\documentclass[12pt,a4paper]{article}
\usepackage{style}
\begin{document}
\maketitle
\medbreak

\begin{abstract}
\noindent We obtain tight bounds for the minimal number of generators of an ideal with bounded-degree generators in a polynomial ring $K[X_1,\dots,X_n],$ as well as a sharp quantification of the maximum possible size of a minimal generating set of bounded degree. Our bounds are sharp for all fields of size greater than the degree. Moreover, we provide explicit constructions reaching the tightness constraints for all fields of characteristic $0,$ and for all sufficiently large fields in the one- and two-variable case. Additionally, we fully solve the one-variable case, and conjecture the asymptotics in the multivariate case, for all finite fields.
\end{abstract}
\section{Introduction}
Bounding the number of low degree generators of polynomial ideals is important for various applications, in particular in recent research in spectral theory \cite{avila}. There are two main formulations of this problem: 
The first is as follows. Suppose that $I\subset K[X_1,X_2,\dots,X_n],$ for some field $K$, is an ideal that is generated by some (possibly infinite) set of polynomials, each of degree at most $d.$ By the Noetherian property, this ideal is necessarily finitely generated. Denote by $\mu(I)$ the minimal possible number of generators of $I.$ The constant $\mu(I)$ has been well-studied, and is the subject of the famed Murthy's conjecture $\mu(I)=\mu(I/I^2)$ for all ideals $I$ in polynomial rings, which has recently been proved over all infinite fields \cite{fasel,mandal}. Our first result bounds $\mu(I),$ for ideals $I$ generated by polynomials of degree at most $d,$ from above.
\begin{theorem}\label{thm-1}
Given an ideal $I=(f_j)_{j\in J}\subset K[X_1,X_2,\dots,X_n]$ such that $\deg\, f_j\leq d$ for all $j\in J,$ we have $\mu(I)\leq\binom{n+d-1}{d},$ with all generating polynomials chosen to all have degree at most $d.$ This bound is sharp, i.e. equality is achieved, for all choices of $K,n,d.$
\end{theorem}
An important variation of this problem, and the main focus of our paper, is the one of finding a finite number of generators \emph{from the original set} of generating polynomials $f_j.$ In other words, we wish to bound the size of a \emph{minimal generating (sub)set} of $\{f_j\}_{j\in J},$ where minimality is with respect to set containment. To the best of our knowledge, this problem has not been studied in the literature, but has been briefly mentioned on public forums \cite{stackexchange}\footnote{The cited Math StackExchange post suggests correct proofs for the sharpness result of \Cref{thm-1} and the upper bound given by \Cref{nullstel-thm}. However, the author is not aware of any results in the literature giving the other sides of the bounds.}. 

From \Cref{thm-1} (and as noted in \cite{stackexchange}) we certainly obtain a lower bound of $\binom{n+d-1}{d}.$ Luckily, the upper bound is not too much higher.
\begin{theorem}\label{nullstel-thm}
Every minimal generating set $S\subset K[X_1,\dots,X_n]$ such that $\deg f\leq d$ for all $f\in S$ satisfies $\abs{S}\leq\binom{n+d}{d}.$ This bound is sharp, i.e. equality is achieved, for all $K,n,d$ such that $\abs{K}>d$ (in particular, if $K$ is infinite).
\end{theorem}
The above two results fully solve the problem for all infinite, and sufficiently large finite, fields. However, when $K=\F_q$ is a fixed finite field, and as $d$ increases, a rich theory is exposed: Our paper only touches the tip of the iceberg. \begin{proposition}\label{n-equal-1}
Any minimal generating set $S\subset\F_q[X]$ such that $\deg f\leq d$ for all $f\in S$ satisfies $\abs{S}\leq\lambda_q(d)=\frac{dq^{O(1)}}{\log d}.$ Equality is achieved for all $q$ by an explicit construction.
\end{proposition}
For $n>1,$ we propose (and justify) the following conjecture generalizing \Cref{n-equal-1}:
\begin{conjecture}\label{n-greater-1}
    Every minimal generating set $S\subset\F_q[X_1,\dots,X_n]$ such that $\deg f\leq d$ for all $f\in S$ satisfies $\abs{S}\leq\Theta_{n,q}\left(d/\log d\right)^n,$ and equality is achieved for all $n,q.$
\end{conjecture}
For fixed $n,q,$ these asymptotics are in between the lower bound given by \Cref{thm-1}, of size $\binom{n+d-1}{d}=\Theta_n(d^{n-1}),$ and the upper bound $\binom{n+d}{d}=\Theta_n(d^n)$ given by \Cref{nullstel-thm}. The author has argumentation nearing (but missing a key idea for) a proof of the equality case of \Cref{n-greater-1}. 

In the next section, we prove \Cref{thm-1,nullstel-thm}. Following this, in \Cref{finite-field}, we prove the asymptotics of \Cref{n-equal-1} and develop arguments towards proving \Cref{n-greater-1}. In the final \Cref{explicit-constructions}, we provide explicit constructions meeting the sharp bound of \Cref{nullstel-thm} in all cases when $K$ has characteristic $0,$ and for all sufficiently large fields when $n\leq2.$
\section{Proofs of Main Theorems}
The proof of \Cref{thm-1} consists of a simple telescoping argument. Precisely, the existence of lower-degree generators of any given ideal allows for a decreased number of higher-degree generators.
\begin{proof}[Proof of \Cref{thm-1}]
Fix an ideal $I$ generated by polynomials of degree at most $d.$ Let $V_{\leq k}$ denote the subspace of $I$ consisting of polynomials of degree at most $k.$ Let $V_k$ be the set of all degree-$k$ parts of elements of $V_{\leq k},$ which is also a $K$-linear space; note however that it is not necessarily a subspace of $V_{\leq k}.$ Let $c_k=\dim_KV_k.$ We will prove by induction that $( V_{\leq k})$ can be generated by $c_k$-many polynomials of degree at most $k$ for all $k.$ Since $I=( V_{\leq d})$ and $c_k\leq\binom{n+k-1}{k}$ (as a subspace of homogeneous degree-$k$ polynomials), this would imply our conclusion.

The base case of $k=0$ is trivial; either $c_0=0$ in which case $V_{\leq k}=\{0\}$ or $c_0=1$ in which case $V_{\leq k}=K$ is generated by $1.$ Now, suppose the hypothesis holds true for a fixed $k.$ Then fix $c_k$-many generators $f_i$ of degree at most $k$ of $( V_{\leq k}).$ Now, we know that for every $f\in V_k$ there exists a polynomial combination of the $f_i$ with degree-$k$ part equal to $f.$ In particular, multiplying by $X_1,$ we conclude that for every $f\in X_1V_k$ there exists a polynomial combination of $f_i$ with degree-$(k+1)$ part equal to $f.$ We have that $X_1V_k\subset V_{k+1}$ has dimension $c_k,$ and so by standard linear algebra we can find $c_{k+1}-c_k=d_{k+1}$-many $\a_j\in V_{k+1}$ that, together with $X_1V_k,$ generate $V_{k+1}$ as a $K$-vector space. For each $\a_j,$ choose $g_j\in V_{\leq(k+1)}$ such that the degree-$(k+1)$ part of $g_j$ is precisely $\a_j.$ Then from here it is simple to see that in fact $g_j,$ together with $f_i,$ generate $( V_{\leq(k+1)}).$ There are $d_{k+1}+c_k=c_{k+1}$-many total of them, so the inductive step holds, and thus a solution with at most $\binom{n+d-1}{d}$ generators indeed exists. In fact, we have shown that a solution with at most $c_d$-many generators exists.

Sharpness is much easier to see; indeed, as suggested in \cite{stackexchange}, consider the ideal $I=(X_1,X_2,\dots,X_n)^d$ generated by all monomials of degree exactly $d$; note that there are precisely $\binom{n+d-1}{d}$-many of these. Suppose that $I$ can be generated by monomials $f_1,f_2,\dots,f_k$ of degree at most $d.$ Observe that since $f_i\in I,$ and $f_i$ has degree at most $d,$ every single monomial with nonzero coefficient in $f_i$ must have degree precisely $d$ (the degree-$0$ monomial $1$ cannot be present, as only the degree-$d$ monomials are generators of $I$).

Now, if $\sum_ig_if_i=m$ for some monomial $m$ of degree $d,$ then comparing degree-$d$ parts, and letting $c_i=g_i(\mathbf{0})$ be the constant coefficients, we see that $\sum_ic_if_i=m$ and thus a $K$-linear combination of $f_i$ equals $m.$ Hence, in fact $f_i$ generate the $K$-vector space of homogeneous degree-$d$ polynomials. This space has dimension $\binom{n+d-1}{d},$ hence we conclude that $k\geq\binom{n+d-1}{d},$ as desired.
\end{proof}
The statement of \Cref{nullstel-thm} (and, similarly, \Cref{n-equal-1}) can be rewritten as follows: Any generating set of degree-at-most-$d$ polynomials of an ideal has a generating subset of size at most $\binom{n+d}{d}.$ Or, more formally: Given an ideal $I=(f_j)_{j\in J}\subset K[X_1,X_2,\dots,X_n],$ there exists $J_*\subset J$ such that $\{f_j:j\in J_*\}$ also generate $I,$ and $\abs{J_*}\leq\binom{n+d}{d}.$ Our
proof uses the idea of evaluations of polynomials at points. Precisely, if for each generating polynomial $f_j$ there exists a point $P_j$ at which $f_j$ is nonzero and all other generators evaluate to $0,$ then $f_j$ is an `essential' generator that cannot be removed.
\begin{proof}[Proof of \Cref{nullstel-thm}]
As suggested in \cite{stackexchange}, the upper bound can be proved as follows. Say that $I=( f_i)_{i\in I}$ is our ideal. Observe $\{f_i\}$ generate a linear subspace of $V=K[X_1,X_2,\dots,X_n]_{\leq d},$ the space of all $n$-variable polynomials of degree at most $d.$ This subspace has dimension $D\leq\dim V=\binom{n+d}{d}$ and thus we can choose a basis of size $D$ for it among the $f_i,$ say we choose $f_1,f_2,\dots,f_D.$ Then each $f_i$ is a $K$-linear combination of $f_1,f_2,\dots,f_D$ and thus in particular lies inside the ideal spanned by them. So, we in fact conclude that $I=( f_1,f_2,\dots,f_D)$ as desired.

For sharpness, the key idea is to consider the \emph{evaluation maps} $\phi_P$ at points $P\in K^n,$ which are linear maps (and in fact also ring homomorphisms) from $K[X_1,\dots,X_n]$ to $K.$ One standard method of proving that a certain set $J$ indexing generators $f_j$ of an ideal $I$ is minimal is as follows: Suppose for each $j\in J$ we designate a point $P_j\in K^n.$ Now, suppose that $f_k(P_j)=0$ if $j\neq k,$ but $f_j(P_j)\neq0$ for all $j\in J.$ Then, observe that no proper subset $J'\subset J$ has the property that $\{f_j:j\in J'\}$ generate all of $I$: Indeed, fix such a set $J',$ and fix $j_*\in J-J'.$ Then $( f_j:j\in J')\subset\Ker\phi_{P_{j_*}}$ (note that the latter is often instead written as the maximal ideal $\mm_{P_{j_*}}$) as $f_j(P_{j_*})=0$ for all $j\in J'.$ On the other hand $f_{j_*}\not\in\Ker\phi_{P_{j_*}},$ and thus $I\not\subset\Ker\phi_{P_{j_*}}.$ Consequently $\{f_j:j\in J'\}$ do not generate $I.$

Hence, it remains to choose $\binom{n+d}{d}$-many points $P_j,$ and consequently for each $j$ construct a polynomial $f_j\in K[X_1,\dots,X_n]$ of degree at most $d$ such that $f_j(P_j)\neq0$ but $f_j(P_k)=0$ for all $k\neq j.$ Here, the combinatorial Nullstellensatz proves quite useful. Indeed, since $\phi_P$ is linear, it can in fact be dually presented via a vector $v_P,$ with the action $\phi_P(f)$ simply being the inner product of $v_P$ with the vector of coefficients of $f.$ Explicitly, if we write the coefficients $$
c_f:=\begin{pmatrix}
c_1&c_{X_1}&c_{X_2}&\hdots&c_{X_n}&c_{X_1^2}&c_{X_1X_2}&\hdots&c_{X_n^d}
\end{pmatrix}
$$
then if $P=(p_1,p_2,\dots,p_n)$ then we have $$
v_P=\left(\begin{matrix}1&p_1&p_2&\cdots&p_n&p_1^2&p_1p_2&p_2^2&\cdots&p_n^d\end{matrix}\right.\kern+0.4em)^\top.
$$
Hence, given a set of points $P_1,P_2,\dots,P_{\binom{n+d}{d}},$ we in fact have the identity $$
\begin{pmatrix}
    \phi_{P_1}(f)\\\phi_{P_2}(f)\\\vdots\\\phi_{P_{\binom{n+d}{d}}}(f)
\end{pmatrix}
=
\begin{pmatrix}
    v_{P_1}^\top\\v_{P_2}^\top\\\vdots\\v_{P_{\binom{n+d}{d}}}^\top
\end{pmatrix}c_f.
$$
Let $A=A_{P_1,P_2,\dots,P_{\binom{n+d}{d}}}:=\begin{pmatrix}
    v_{P_1}^\top\\v_{P_2}^\top\\\vdots\\v_{P_{\binom{n+d}{d}}}^\top
\end{pmatrix}$, a so-called \emph{multivariate Vandermonde matrix} (useful in the study of Reed-Muller codes \cite{reedmuller}). Now, we wish to find $f_i,\,1\leq i\leq\binom{n+d}{d},$ such that $Ac_{f_i}=e_i,$ the standard basis vector with entry $1$ in the $i$-th coordinate. Evidently, these exist for all $i$ if and only if the matrix $A$ is invertible. Its determinant is a nonzero (in fact, every permutation gives a distinct monomial) degree-$\sum\limits_{i=0}^di\binom{n+i-1}{i}=\binom{n+d}{d-1}$ polynomial in the coordinates $p_{ij}$ of the points $P_i.$ In fact note that the maximum possible degree of any one of the $p_{ij}$ in any monomial is $d,$ since the term $p_{ij}$ appears only in the $i$-th row of $A.$ Hence, by the combinatorial Nullstellensatz \cite{alon}, since $\abs{K}>d,$ there exists a choice of $p_{ij}\in K$ for all $1\leq i\leq\binom{n+d}{d},1\leq j\leq n,$ such that $\det\,A\neq0.$ 

Consequently, there exists a choice of (trivially distinct) points $P_1,\,P_2,\,\dots,\,P_{\binom{n+d}{d}}\in K^n$ for which there exist respective polynomials $f_1,\,f_2,\,\dots,\,f_{\binom{n+d}{d}}$ of degree at most $d$ such that $f_i(P_j)=\delta_{ij}=\begin{cases}1&i=j,\\0&i\neq j.\end{cases}$ Then, by our previous argument, the ideal $I=( f_i:1\leq i\leq\binom{n+d}{d})$ requires all $\binom{n+d}{d}$-many generators, as desired.
\end{proof}
\section{Asymptotics for Finite Fields}\label{finite-field}
In \Cref{nullstel-thm}, a sharp bound is achieved for all fields $K$ with $\abs{K}>d.$ What happens if $\abs{K}\leq d$? More specifically, the following question arises: suppose $K=\F_q$ for some fixed prime power $q,$ and $n$ is fixed. What is the bound $m(n,d),$ asymptotically as $d\to\infty$? We restate \Cref{n-equal-1}, which estimates this bound for $n=1$:
\begin{proposition*}
    Any minimal generating set $S\subset\F_q[X]$ such that $\deg f\leq d$ for all $f\in S$ satisfies $\abs{S}\leq\lambda_q(d)=\frac{dq^{O(1)}}{\log d}.$ Equality is achieved for all $q$ by an explicit construction.
\end{proposition*}
We can only approximate $\lambda_q(d)$ by an asymptotic that varies within a constant factor, depending on the exact value of $d.$ Regardless, this theorem proves that the growth of the largest possible minimal generating set of an ideal in $\F_q[X]$ is sublinear.

The reason we focus first on $n=1$ is that for any (even finite) field $K,$ the ring $K[X]$ is a Euclidean domain, allowing for arguments involving divisibility. Indeed, suppose $I\subset K[X]$ is an ideal with a minimal generating set $(f_1,f_2,\dots,f_m).$ Our goal is to lower bound the maximum degree of the $\{f_i\},$ in terms of $m.$ We may assume without loss of generality that $\{f_i\}$ do not all share a common factor, as otherwise we can divide all of them by this factor to obtain a lower-degree example. In particular, by the Euclidean domain property, in fact $$
(f_1,f_2,\dots,f_m)=(\gcd(f_1,f_2,\dots,f_m))=(1)=K[X].
$$
In particular, since $\{f_i\}$ is a minimum generating set, this means that for each $i,$ the set $\{f_1,\dots,f_{i-1},f_{i+1},\dots,f_m\}$ does not generate all of $K[X].$ Again, by the Euclidean domain property, this implies that these $(m-1)$-many polynomials share a common irreducible factor; call it $g_i.$ Observe that $g_i\nmid f_i$ since otherwise $g_i$ divides all of $f_1,f_2,\dots,f_m.$ 

In particular, we easily conclude that $g_i$ are all distinct irreducible polynomials, and for each $i$ the polynomial $f_i$ must therefore be divisible by the product $\prod\limits_{j\neq i}g_j.$ 

In fact, we can set $f_i=\prod\limits_{j\neq i}g_j$ for each $i,$ and this will indeed be a minimal generating set. True to this fact, observe that indeed $g_i$ divides $f_j$ if and only if $i\neq j.$ Hence, if any $f_i$ is removed from the set, then the remaining polynomials generate an ideal contained in $(g_i),$ a proper ideal of $K[X].$ In other words, we conclude that indeed $(f_1,f_2,\dots,f_m)$ is a minimal generating set.

In the case $\abs{K}>d,$ observe that there are $(d+1)$-many distinct elements $\a_0,\a_1,\dots,\a_d\in K,$ and set $g_i=X-\a_i,$ which are distinct irreducibles. Then each $f_i$ is a product of $d$-many linear polynomials, hence of degree $d.$ In other words, we have constructed a minimal generating set of size $d+1=\binom{1+d}{d}$ and degree at most $d,$ as desired. This notably solves the problem for all infinite fields.

However, if $K=\F_q$ is finite, then we may enumerate the irreducible polynomials in $\F_q[X]$ in order of non-decreasing degree as $h_1,h_2,\dots.$ In our above minimal-degree example, let us assume without loss of generality that $g_i$ are ordered by increasing degree, so that necessarily $\deg g_i\geq\deg h_i$ for all $i.$ Then observe that $f_1$ is divisible by the product $\prod\limits_{i=2}^mg_i$ and thus $\deg\,f_1\geq\sum\limits_{i=2}^{m}\deg\,g_i\geq\sum\limits_{i=2}^{m}\deg\,h_i.$ In order to complete the proof of \Cref{n-equal-1}, it thus remains to estimate the growth of $\deg h_m$ in $m.$
\begin{proof}[Proof of \Cref{n-equal-1}]
    By Gauss's formula counting irreducible polynomials of degree $k$ over the finite field $\F_q$ (see a new derivation by Chebolu and Min\'a{\v{c}} \cite{chebolu}): $$
    p_q(k)=\frac{1}{k}\sum\limits_{d\mid k}\mu(d)q^{k/d}=\frac{q^k+O(q^{k/2})}{k}.$$
    where $\mu$ is the M\"obius function. In particular, summing over $k,$ we may conclude that the number of irreducible polynomials of degree \emph{at most} $n$ is \begin{equation}\label{eqn:irreducible-bound}
    P_q(n)=\frac{q^{n+1}+O(q^{n}/n)}{n(q-1)},
    \end{equation}
    Now, by our argument immediately above, observe that if $\{f_1,f_2,\dots,f_m\}$ is a minimal generating set, then $\deg\,f_1\geq\sum\limits_{i=2}^m\deg\,h_i.$ On the other hand, we know that $\deg\,f_1\leq d.$ Let $n$ be the largest integer such that $P_q(n)\leq m.$ Observe by Equation \ref{eqn:irreducible-bound} that in fact $n=\log_qm+\log_q\log_qm+O(1).$ Correspondingly $P_q(n)=mq^{-O(1)}.$ Then, for each $1\leq j\leq n,$ there are exactly $p_q(j)$-many of the $h_i$ which have degree $j$; the remaining $(m-P_q(n))$-many polynomials all have degree $n+1.$ In other words we conclude that (since $h_1$ is missing)
    \begin{align*}d&\geq\sum\limits_{j=1}^njp_q(j)+(n+1)(m-P_q(n))-1\\
    &=\sum\limits_{j=1}^n\frac{q^{j+1}+O(q^{j/2})}{q-1}+m(\log_qm+\log_q\log_qm+O(1))(1-q^{-O(1)})\\
    &=\frac{q^{n+2}}{(q-1)^2}+(1-q^{-O(1)})m\log_qm+O(q^{n/2})=q^{O(1)}m\log m.
    \end{align*}
    Inverting this bound, we obtain that $m\leq\frac{dq^{O(1)}}{\log d}$ as desired.

    For sharpness, observe that we can use the construction we discussed for the case that $\abs{K}>d$ and apply it to the same $m$ lowest-degree irreducible polynomials $h_1,h_2,\dots,h_m.$ Setting $f_i=\prod\limits_{j\neq i}h_j,$ we conclude by the Euclidean property that $\{f_1,f_2,\dots,f_m\}$ is a minimal generating set. By setting $m$ maximal such that $\sum\limits_{j=1}^njp_q(j)+(n+1)(m-P_q(n))-1\leq d$ (so that $\deg f_1\leq d$; it is not hard to see that $\deg f_i\leq\deg f_1$ for all $1\leq i\leq m$), we obtain our desired matching lower bound.
\end{proof}
\noindent The above arguments solve the problem entirely for $n=1.$ However, for $n>1,$ we may not use the Euclidean property anymore. This makes the problem, especially for finite fields, rather more difficult, as minimality of a generating set is no longer equivalent to a criterion of divisibility by irreducible polynomials. 

Yet, the use of irreducible polynomials motivates us to look at points in field extensions. Indeed, the solutions of an irreducible polynomial over $\F_q$ of degree $k$ lie in the field $\F_{q^k}.$ In particular, our problem can in fact still be viewed linearly, but the linear factors lie in field extensions $\F_{q^k}[X],$ and multiply to elements of $\F_q[X].$ We can try to generalize this idea to the case of $n>1,$ which leads us to \Cref{n-greater-1}. We now restate it here:

\begin{conjecture*}
Every minimal generating set $S\subset\F_q[X_1,\dots,X_n]$ such that $\deg f\leq d$ for all $f\in S$ satisfies $\abs{S}\leq\Theta_{n,q}\left(d/\log d\right)^n,$ and equality is achieved for all $n,q.$
\end{conjecture*}

Our heuristics prompting this conjecture use an unproven variation of \Cref{nullstel-thm} which relies on a density assumption. Precisely, let us choose $k=\lfloor\log_qd\rfloor+1,$ so that $q^k>d/k.$ Set $d'=\lfloor d/k\rfloor,$ so $q^k>d'.$ By \Cref{nullstel-thm}, there exist (distinct) points $P_1,\dots,P_{\binom{n+d'}{d'}}\in\F_{q^k}^n$ and respective degree-at-most-$d'$ polynomials $g_1,g_2,\dots,g_{\binom{n+d'}{d'}}\in\F_{q^k}[X_1,\dots,X_n]$ such that $g_i(P_j)=\delta_{ij}$ for all $i,j.$ Our \emph{unproven assumption}, however, will be that this choice can be adjusted to additionally satisfy $g_i(\sigma(P_i))\neq0$ for all $i$ and $\sigma\in\Gal(\F_{q^k}/\F_q).$


Assuming the existence of such polynomials and respective evaluation points, since $\F_{q^k}$ is a Galois extension of $\F_q,$ we can simply define $$
f_i=N_{\F_{q^k}/\F_q}(g_i)=\prod\limits_{\sigma\in\Gal(\F_{q^k}/\F_q)}\sigma(g_i)\in\F_q[X_1,\dots,X_n].
$$
where the automorphisms $\sigma$ are applied to the coefficients. Indeed $f_i\in\F_q[X_1,\dots,X_n]$ since their coefficients are fixed by every automorphism in $\Gal(\F_{q^k}/\F_q),$ and thus lie in the fixed field $\F_q.$ 

Observe that for any $i\neq j,$ we have $f_i(P_j)$ is divisible by $g_i(P_j)=0,$ and thus $f_i(P_j)=0.$ At the same time, for all $i$ we have $$
f_i(P_i)=\prod\limits_{\sigma\in\Gal(\F_{q^k}/\F_q)}\s(g_i)(P_i)=\prod\limits_{\sigma\in\Gal(\F_{q^k}/\F_q)}\s^{-1}(g_i(\s(P_i)))\neq0.
$$
Each $P_i$ gives a maximal ideal $\mm_{P_i}=\{f\in\F_q[X_1,\dots,X_n]:f(P_i)=0\}.$ We then have $f_i\not\in\mm_{P_i}$ but $f_j\in\mm_{P_i}$ for all $j\neq i.$ Hence by the same argument as the one given in \Cref{nullstel-thm} we may conclude that $\{f_i\}$ is indeed a minimal generating set. We have $\deg f_i=k\deg g_i\leq d$ and there are $$
\binom{n+d'}{d'}=\Theta_n(d')^n=\Theta_{n,q}(d/\log d)^n
$$
polynomials $f_i,$ precisely the number specified by the bound of \Cref{n-greater-1}. Thus establishing the sharpness aspect of this conjecture amounts to proving the density assumption.

\section{Explicit Constructions}\label{explicit-constructions}

    In the proof of \Cref{nullstel-thm}, the existence of ideals proving the sharpness of the bound is entirely non-explicit, and instead relies on the invertibility of a multivariate Vandermonde matrix. In the case of an infinite field $K,$ given a generic choice of evaluation points $P_i,$ it is most probable that polynomials $f_i$ of degree at most $d$ satisfying $f_i(P_j)=\delta_{ij}$ can be found. Indeed, the set of points $P_i$ for which such interpolating polynomials do not exist lies on a union of hypersurfaces of total degree at most $\binom{n+d}{d-1}.$ However, even so, this is not a guarantee, and moreover finding the polynomials $f_i$ after choosing points $P_j$ is quite tedious (the beautiful Vandermonde polynomial factorization formula does not generalize to the multivariate case; see \cite{brown} for a recent development decomposing the determinant into a sum of fully-factorized terms). This merits looking instead for \emph{explicit} constructions of ideals meeting the sharp bound. And indeed, this study opens up an interesting theory.
\begin{proposition}
    If $K$ has characteristic $0,$ then an ideal $I\subset K[X_1,\dots,X_n]$ with minimal generating set $\{f_1,\dots,f_{\binom{n+d}{d}}\}$ can be explicitly constructed by considering evaluation points $(d_1,d_2,\dots,d_n)\in\N_0^n$ satisfying $d_1+\cdots+d_n\leq d.$
\end{proposition}
\noindent For example, in the case $n=2$ and $d=2,$ the proof of the theorem will give us that the ideal $( XY,X(X-1),Y(Y-1),X(X+Y-2),Y(X+Y-2),(X+Y-1)(X+Y-2))$ is minimally generated by these generators. 
\begin{remark}
    In general, the $\binom{n+d}{d}$-many evaluation points form an $n$-dimensional discrete simplex of depth $d+1,$ and the polynomials constructed all have $d$-many linear factors; for each factor, the condition that it is nonzero `peels off' one layer of the simplex.
\end{remark}
\begin{proof}
Fix a point $P=(d_1,d_2,\dots,d_n)$ of the form specified above. Suppose $f=\sum\limits_{i=1}^nd_i,$ and define the (indeterminate) $X=\sum\limits_{i=1}^nX_i.$ Then, consider the polynomial $$
f_P:=\left(\prod\limits_{i=1}^n\prod\limits_{j=0}^{d_i-1}(X_i-j)\right)\prod\limits_{i=f+1}^d(X-i).
$$
Its degree is $\sum\limits_{i=1}^nd_i+(d-f)=d$ by definition of $f.$ We claim that $f_P(P)\neq0,$ whereas $f_P(Q)=0$ for any other $Q$ among the evaluation points. Indeed, we have $$
f_P(P)=\left(\prod\limits_{i=1}^n\prod\limits_{j=0}^{d_i-1}(d_i-j)\right)\prod\limits_{i=f+1}^d(f-i)=(-1)^{d-f}(d-f)!\left(\prod\limits_{i=1}^nd_i!\right)\neq0.
$$
On the other hand, fix any other point $Q=(e_1,e_2,\dots,e_n)\neq P.$ Let $e=\sum\limits_{i=1}^ne_i.$ Suppose first that there exists $i_*$ such that $e_{i_*}<d_{i_*}.$ In particular, there exists $0\leq j<d_{i_*}$ such that $e_{i_*}-j=0.$ Then observe that $\prod\limits_{j=0}^{d_{i_*}-1}(e_{i_*}-j)=0,$ and thus $$
f_P(Q)=\left(\prod\limits_{i=1}^n\prod\limits_{j=0}^{d_i-1}(e_i-j)\right)\prod\limits_{i=f+1}^d(e-i)=0.
$$
If $e_i\geq d_i$ for all $i,$ then since $Q\neq P,$ this means that the inequality is strict for some $i,$ and thus in fact $e=\sum\limits_{i=1}^ne_i>\sum\limits_{i=1}^nd_i=f.$ But since $Q$ is one of our evaluation points, we also have $e\leq d.$ Hence, there exists $f<i\leq d$ such that $i=e,$ or equivalently $e-i=0.$ Consequently, the second factor of the above product $\prod\limits_{i=f+1}^d(e-i)=0$ and thus once again $f_P(Q)=0.$ In either case $f_P(Q)=0$ and so we have our $\binom{n+d}{d}$-many polynomials, as desired.
\end{proof}
In positive characteristic $p,$ the picture is not as simple. Indeed, the above example collapses as $X_i-p=X_i.$ However, for $n\leq 2,$ we can still construct explicit, interesting examples. The author has reason to believe that his `pyramidal' method does not work beyond $n=2,$ but there are other approaches which may yield nice results.  

The case $n=1$ has been solved in the beginning of the previous section, save for a natural choice of elements $\a_i\in K$ for fields $K$ with $\abs{K}>d.$ For those of characteristic $0,$ we can take $\a_i=i.$ For fields with a transcendental element $q,$ we can consider the $q$-analog $\a_i=\frac{1-q^i}{1-q}.$ This latter choice can be applied to fields $K$ algebraic over $\F_p$ as well: all we need is an element $q$ such that all powers $1,q,q^2,\dots,q^d$ are distinct. In other words, $q\in K^*$ is an element of multiplicative order strictly greater than $d$; such an element exists for any $K$ with $\abs{K}>d+1$ because any finite subgroup of $K^*$ is cyclic. (If $\abs{K}=d+1,$ then we can simply choose all $(d+1)$-many elements of $K,$ in an unspecified order).

For $n=2,$ however, for infinite (and sufficiently large finite) fields, we can write down an explicit generating set matching the upper bound $\binom{n+d}{d}$ as follows.
\begin{theorem}
    For any two nonnegative integers $i,j$ satisfying $i+j\leq d,$ define $P_{ij}=\left(\frac{x^i-y^i}{x^i},\frac{x^j-y^j}{y^j}\right)\in K(x,y)^2.$ If a pair $(x,y)\in(K^*)^2$ can be chosen such that $x^n\neq y^n$ for $1\leq n\leq d,$ then 
    an ideal $I\subset K[X,Y]$ with minimal generating set $\{f_1,\dots,f_{\binom{n+2}{2}}\}$ can be explicitly constructed by considering evaluation points $P_{ij}.$
\end{theorem}
\noindent Again, observe that we are arranging our points $P_{ij}$ in a triangle, and each linear factor in our constructed polynomials will `peel off' an edge of the triangle. This theorem gives us an explicit construction of an ideal $I$ for any field $K$ satisfying $\abs{K}>d+1,$ as then $(x,y)$ can be chosen as $(1,\zeta),$ where $\zeta$ is any element of the multiplicative group $K^*$ with order at least $d+1.$ Observe also that this bound on the size of $K$ nearly meets our non-constructive argument relayed in Theorem \ref{nullstel-thm} (differing only by $1$).
\begin{proof}
    We first prove several linear relations satisfied by the points $P_{ij}.$ Indeed, first observe that for a fixed $i,$ the points $P_{ij}$ for varying $j$ all satisfy $X=\frac{x^i-y^i}{x^{i}}.$ Similarly, for a fixed $j,$ the points $P_{ij}$ over varying $i$ all satisfy $Y=\frac{x^j-y^j}{y^{j}}.$ Finally, if $i+j=k$ is fixed, then we see that $x^kX+y^kY=x^j(x^i-y^i)+y^i(x^j-y^j)=x^{i+j}-y^{i+j}=x^k-y^k.$

    Hence, analogously to the characteristic-zero situation, for each point $P_{ij},$ consider the degree-$d$ polynomial written as the product $$
    f_{ij}(X,Y)=\prod\limits_{i'=0}^{i-1}\left(x^{i'}X-(x^{i'}-y^{i'})\right)\prod\limits_{j'=0}^{j-1}\left(y^{j'}Y-(x^{j'}-y^{j'})\right)\prod\limits_{k=i+j+1}^d(x^kX+y^kY-(x^k-y^k)).
    $$
    We observe that \begin{align*}
    f_{ij}(P_{ij})&=\prod\limits_{i'=0}^{i-1}\left(x^{i'-i}(x^i-y^i)-(x^{i'}-y^{i'})\right)\prod\limits_{j'=0}^{j-1}\left(y^{j'-j}(x^j-y^j)-(x^{j'}-y^{j'})\right)\\
    &\,\,\,\prod\limits_{k=i+j+1}^d\left(x^{k-i}(x^i-y^i)+y^{k-j}(x^j-y^j)-(x^k-y^k)\right)\\
    &=\prod\limits_{i'=0}^{i-1}\left(y^{i'}-x^{i'-i}y^i\right)\prod\limits_{j'=0}^{j-1}\left(y^{j'-j}x^j-x^{j'}\right)\prod\limits_{k=i+j+1}^d\left(y^{k-j}x^j-x^{k-i}y^i\right)
    \end{align*}
    In particular, as long as $x^{i'-i}\neq y^{i'-i},\,x^{j'-j}\neq y^{j'-j},\,x^{k-i-j}\neq y^{k-i-j},$ for all choices of $0\leq i'<i,\,0\leq j'<j$ and $i+j<k\leq d,$ all the factors are consequently nonzero and thus $f_{ij}(P_{ij})\neq0.$ But this set of conditions is in fact equivalent to $x^n\neq y^n$ for $1\leq n\leq d,$ as mentioned in the statement of our theorem.

    On the other hand, suppose $(\Tilde{i},\Tilde{j})\neq(i,j)$; first assume that one of the coordinates of the first pair is less than the corresponding coordinate of the second pair, without loss of generality $\Tilde{i}<i.$ Now expanding and simplifying, we have $$
    f_{ij}(P_{\Tilde{i}\Tilde{j}})=\prod\limits_{i'=0}^{i-1}\left(y^{i'}-x^{i'-\Tilde{i}}y^{\Tilde{i}}\right)\prod\limits_{j'=0}^{j-1}\left(y^{j'-\Tilde{j}}x^{\Tilde{j}}-x^{j'}\right)\prod\limits_{k=i+j+1}^d\left(y^{k-\Tilde{j}}x^{\Tilde{j}}-x^{k-\Tilde{i}}y^{\Tilde{i}}\right)=0,
    $$
    since in the first product, there exists a choice of $0\leq i'<i$ such that $i'=\Tilde{i},$ and thus one of the factors is $y^{\Tilde{i}}-x^{\Tilde{i}-\Tilde{i}}y^{\Tilde{i}}=y^{\Tilde{i}}-y^{\Tilde{i}}=0.$ 

    Otherwise, we have $\Tilde{i}\geq i,\Tilde{j}\geq j,$ but one of the inequalities is strict. In particular we have $\Tilde{i}+\Tilde{j}>i+j.$ Consequently there exists a choice of $i+j<k\leq d$ such that $k=\Tilde{i}+\Tilde{j},$ and thus again, in the above expression, but now in the third product, one of the factors is $y^{\Tilde{i}+\Tilde{j}-\Tilde{j}}x^{\Tilde{j}}-x^{\Tilde{i}+\Tilde{j}-\Tilde{i}}y^{\Tilde{i}}=y^{\Tilde{i}}x^{\Tilde{j}}-x^{\Tilde{j}}y^{\Tilde{i}}=0,$ so the entire expression is $0$ in this case as well.

    Therefore, assuming the condition on $x$ and $y,$ we conclude that $f_{ij}(P_{\Tilde{i}\Tilde{j}})=0$ if and only if $i=\Tilde{i}$ and $j=\Tilde{j}.$ Notably, since we have constructed polynomials $f_{ij}$ for each pair $(i,j),$ this means that $P_{ij}$ are all distinguished from each other by the action of these polynomials on them. Hence, in fact $P_{ij}$ are $\binom{d+2}{2}$-many distinct points, and the ideal $(\{f_{\Tilde{i}\Tilde{j}}:(\Tilde{i},\Tilde{j})\neq(i,j)\})$ evaluates to $0$ at $P_{ij},$ whereas the ideal $I=(\{f_{\Tilde{i}\Tilde{j}}\})$ does not (since it contains $f_{ij}$ as a generator). Hence, every generator of the ideal $I$ is necessary, and thus $\{f_{ij}\}$ is a minimal generating set of size $\binom{d+2}{2}$ and degree bounded by $d,$ as desired.
\end{proof}
\noindent The simplex idea unfortunately does not extend beyond $n=2,$ which we now explain. Indeed, for $n=3,$ suppose after an affine shift that our initial points are $$
P_{000}=(0,0,0),P_{100}=(x,0,0),P_{010}=(0,y,0),P_{001}=(0,0,z).
$$
Then, necessarily (in order for points to lie on intersections of specified planes), we must also have the identities $$
P_{110}=(x,y,0),P_{101}=(x,0,z),P_{011}=(0,y,z).
$$
These points should all lie on one plane, which should also contain $P_{200},P_{020},P_{002}.$ The plane has equation $aX+bY+cZ=d,$ and solving for $a,b,c,d,$ we obtain that $ax+by=ax+cz=by+cz,$ or equivalently $ax=by=cz,$ so $a=\frac{t}{x},b=\frac{t}{y},c=\frac{t}{z}$ for some $t\in K.$ After rescaling, we obtain that the unique plane on which these three points all lie is $\frac{X}{x}+\frac{Y}{y}+\frac{Z}{z}=2.$ In particular, from here we conclude that $$
P_{200}=(2x,0,0),P_{020}=(0,2y,0),P_{002}=(0,0,2z).
$$
Continuing with this argument inductively, we are forced to set $P_{ijk}=(ix,jy,kz)$ for each ordered triple $(i,j,k).$ However, $K$ has some characteristic $p.$ In particular $P_{p00}=(px,0,0)=(0,0,0)=P_{000}.$ This is not suitable since our evaluation points are required to be distinct.

However, note that for the two distinct problems studied (\Cref{thm-1,nullstel-thm}), the quality of ideals used is rather different. The sharpness bound for the minimum number of generators (not necessarily from the original ones) of an ideal is achieved for the ideal of the form $( X_1,X_2,\dots,X_n)^d,$ whose radical has a single point (the origin) at which it evaluates to $0.$ Correspondingly, proving that $\binom{n+d-1}{d}$ is the correct minimal number of generators of this ideal is a problem of significantly distinct quality than quantifying the minimal number of generators among $\binom{n+d}{d}$ specified ones; indeed, the explicit ideal $I$ 
that we have constructed in each case for the latter problem is in fact generated by one element, the constant $1$ (indeed, if there is a minimal set of $\binom{n+d}{d}$ generators, then they must be linearly independent and thus span all of $K[X_1,\dots,X_n]_{\leq d},$ which includes $1$). Our argument inherently relies on the fact that in fact $V(I)$ is an (empty) variety which is the intersection of hypersurfaces $V(f_j)$ corresponding to each generator $f_j$ of $I$; since for any $i\in J$ the intersection of $V(f_j)$ over all $j\neq i$ is nonempty, this implies that all sets are required to be included in order to have overall empty intersection. Therefore, for $n>2$ and in characteristic $p,$ there is still hope to come up with an explicit description of a minimal generating set meeting the sharp upper bound, if we consider ideas beyond, or perhaps combining, the two described previously.
\section*{Acknowledgements}
We are grateful to Artur Avila for the question that led to this work. We also thank Melanie Matchett Wood and David Eisenbud for helpful comments on the earlier version of this paper that have led to a significant improvement.

\bibliographystyle{roman}

\begin{thebibliography}{1}
    \bibitem{alon} Noga Alon (1999). \emph{Combinatorial Nullstellensatz}. Combinatorics, Probability and Computing. 8 (1–2): 7--29.
    \bibitem{avila} Artur Avila, David Damanik (2025). \emph{Pure Point Spectrum is Generic}. Preprint.
    \bibitem{brown} Francis Brown (2025). \emph{Multivariable Vandermonde determinants, amalgams of matrices and Specht modules}. Journal of Algebra. 678: 253--278. 
    \bibitem{chebolu} Sunil K.~Chebolu, J\'an Min\'a{\v{c}} (2011). \emph{Counting irreducible polynomials over finite fields using the inclusion-exclusion principle}. Mathematics Magazine, 84 (5): 369--371.
    \bibitem{fasel} Jean Fasel (2016). \emph{On the number of generators of ideals in polynomial rings}. Annals of Mathematics, 184: 315--331.
    \bibitem{mandal} Satya Mandal (2016). \emph{On the complete intersection conjecture of Murthy}. Journal of Algebra. 458: 156--170. 
    \bibitem{reedmuller} Irving S.~Reed (1954). \emph{A class of multiple-error-correcting codes and the decoding scheme}. Transactions of the IRE Professional Group on Information Theory. 4 (4): 38--49.
    \bibitem{stackexchange} Mees de Vries (2023). \emph{Ideal of $F[X_1,\dots,X_n]$ generated by polynomials of bounded degree -- small generating subset?} Math StackExchange, Aug 1, 2023 (\href{https://math.stackexchange.com/questions/4745862/ideal-of-fx-1-dotsc-x-n-generated-by-polynomials-of-bounded-degree-small}{link}).
\end{thebibliography}

\end{document}